\documentclass[11pt, letter]{amsart}

\usepackage[utf8]{inputenc}
\usepackage{amsmath, amssymb, amsthm}
\usepackage{mathtools}
\usepackage{stmaryrd}
\usepackage{enumitem}
\usepackage{graphicx}
\usepackage[usenames,dvipsnames]{xcolor}

\usepackage[backref=page]{hyperref}
\hypersetup{
  colorlinks   = true,
  urlcolor     = teal,
  linkcolor    = teal,
  citecolor    = teal
}
\usepackage[capitalize]{cleveref}
\usepackage{xcolor}

\newtheorem{theorem}{Theorem}[section]
\newtheorem{proposition}[theorem]{Proposition}
\newtheorem{lemma}[theorem]{Lemma}
\newtheorem{corollary}[theorem]{Corollary}
\theoremstyle{definition}
\newtheorem{definition}[theorem]{Definition}
\newtheorem{remark}[theorem]{Remark}

\usepackage[usenames,dvipsnames]{xcolor}

\usepackage[color=Aquamarine,textsize=footnotesize]{todonotes}

\newcommand{\R}{{\bf R}}
\newcommand{\Z}{{\bf Z}}
\newcommand{\dist}{\mathsf{d}}
\newcommand{\Hn}{\mathcal{H}^N}
\newcommand{\til}{\widetilde}
\newcommand{\cur}[1]{\llbracket #1 \rrbracket}
\newcommand{\bI}{\mathbf{I}}
\newcommand{\cI}{\mathcal{I}}
\newcommand{\bM}{\mathbf{M}}
\newcommand{\mres}{%
  \,\raisebox{-.127ex}{\reflectbox{\rotatebox[origin=br]{-90}{$\lnot$}}}\,}
\DeclareMathOperator{\RCD}{\mathsf{RCD}}

\DeclareMathOperator{\CAT}{\mathsf{CAT}}
\DeclareMathOperator{\Jac}{Jac}
\DeclareMathOperator{\spt}{supp}
\DeclareMathOperator{\bary}{bar}
\DeclareMathOperator{\ind}{ind}
\DeclareMathOperator{\Lip}{Lip}
\DeclareMathOperator{\Tan}{Tan}

\begin{document}

\title[Entropy-degree rigidity for non-collapsed $\RCD$ spaces]{Entropy-degree rigidity \\ for non-collapsed $\RCD$ spaces}
\author[P. Su\'arez-Serrato]{P. Su\'arez-Serrato}
\date{\today}

\address{University of California Santa Barbara, Anisq'Oyo - Isla Vista, California}
\address{{\it On leave from}: Instituto de Matem\'aticas, Universidad Nacional Aut\'onoma de M\'exico UNAM}
\address{M\'exico Tenochtitlan}

\begin{abstract}
We prove the entropy-degree theorem for continuous maps from a closed orientable non-collapsed $\RCD(K,N)$ space $X$ to a closed negatively curved locally symmetric $N$-manifold $Y$,  with $N \geq 3$.
In the equality case, under a fundamental group index condition, we show the space is isometric to a smooth locally symmetric Riemannian manifold, and the map is homotopic to a Riemannian covering.
\end{abstract}

\maketitle

\tableofcontents

\section{Introduction}

Let $(Z,\dist,\mathfrak{m})$ be a metric measure space.
Define the \emph{volume entropy} $h(Z)$ as (cf.\ \cite{M, BCG-Samb, Reviron})
\[
h(Z)=\limsup_{R\to\infty}\frac{\log \mathfrak{m}(B(x,R))}{R},
\]
where $B(x,R)$ is the ball of radius $R$ in $Z$ centered at $x$ in $Z$.
Here, $h(Z)$ denotes the exponential volume growth rate of $Z$ itself, and the traditional volume entropy of a compact space is $h(\til{Z})$, computed on its universal cover (cf \cite{RCDbary, SuarezSerrato2026}).
Following a strategy introduced by Sambusetti \cite{Sambusetti-1999}, we use the entropy of intermediate covering spaces.

Given a continuous map $f \colon X \to Y$ between topological spaces $X$ and $Y$, let $f_* \colon \pi_1(X) \to \pi_1(Y)$ be the induced homomorphism, and let $\overline{X}$ be the cover of $X$ associated with the subgroup $\ker f_* < \pi_1(X)$.
Then $\Gamma := \pi_1(X)/\ker f_*$ acts on $\overline{X}$ by deck transformations, and $\overline{X}$ is the smallest cover of $X$ to which $f$ lifts as a map $\overline{f} \colon \overline{X} \to \til{Y}$.
The class of $\RCD(K,N)$ spaces consists of metric measure spaces with synthetic Ricci curvature bounded below by $K$ and dimension bounded above by $N$.
It originates in the curvature-dimension condition of Sturm \cite{Sturm-2006-I, Sturm-2006-II} and Lott--Villani \cite{Lott-Villani-2009}, and in the Riemannian refinement introduced by Ambrosio, Gigli, and Savar\'e \cite{AGS-2014-Duke} and Gigli \cite{Gigli-2015-Memoirs}.
For $\RCD$ spaces, the existence of universal covers was shown by Mondino--Wei \cite{Mondino-Wei-2019}.
Moreover,  Wang \cite{Wang-2024} proved that $\RCD^*(K,N)$ spaces are semi-locally simply connected, and hence $\overline{X}$ is defined by classical covering space theory.

De Philippis and Gigli \cite{DPG} defined the \emph{non-collapsed} $\RCD(K,N)$ spaces, as those for which $N$ is a natural number and the reference measure is the Hausdorff measure $\Hn$.
This subclass contains the non-collapsed Ricci limit spaces studied by Cheeger and Colding \cite{Cheeger-Colding-1997}, and, by the work of Petrunin \cite{petrunin2011} and Zhang--Zhu \cite{zhang2010ricci}, the $N$-dimensional Alexandrov spaces with curvature bounded below, equipped with $\Hn$.
The structure theory of Mondino--Naber \cite{MondinoNaber19}, Kell--Mondino \cite{KellMondino}, and Gigli--Pasqualetto \cite{GigliPasqualetto:CAG} shows that these spaces are $N$-rectifiable with Euclidean tangents almost everywhere. 
The boundary theory of De Philippis--Gigli \cite{DPG}, Kapovitch--Mondino \cite{KapMonGT}, and Bru\`e--Naber--Semola \cite{BrueNaberSemola20} describes their codimension-one singularities.

We say that a non-collapsed $\RCD(K,N)$ space is \emph{closed} if it is compact and its boundary  is empty \cite{DPG}.

Our main result is the following.

\begin{theorem}\label{thm:rcd-entropy-degree}
Let $N \geq 3$ be a natural number, $K$ in $\R$, and $(X,\dist,\Hn)$ a closed non-collapsed $\RCD(K,N)$ space whose topological manifold part $\mathcal{M}(X)$ is orientable.
Let $Y$ be a closed orientable negatively curved locally symmetric space of dimension $N$.
Then, for any continuous map $f \colon X \to Y$,
\begin{equation}\label{eq:rcd-entropy-degree-ineq}
    h(\overline{X})^N\, \Hn(X) \;\geq\; |\deg(f)|\, h(\til{Y})^N\, \Hn(Y),
\end{equation}
where $\deg(f)$ is the degree of $f$ given in \Cref{def:deg-continuous}.

Moreover, if equality holds in inequality \eqref{eq:rcd-entropy-degree-ineq} and $|\deg(f)| = [\pi_1(Y) : f_*\pi_1(X)]$, then, after a homothetic rescaling of the metric on $X$, the space $X$ is isometric to a smooth locally symmetric Riemannian manifold, and $f$ is homotopic to a Riemannian covering map of degree $|\deg(f)|$.
\end{theorem}

The degree in \Cref{thm:rcd-entropy-degree} is defined through local degrees of $f$ on the topological manifold part $\mathcal{M}(X)$ of $X$ (see \Cref{def:deg-Lipschitz}).
It agrees with the classical degree when $X$ is a closed oriented topological manifold.
As in the smooth case, it is a homotopy invariant.
Moreover, by \Cref{prop:degree-equivalence}, it is the multiplicity of the push-forward of the orientation current of $X$ \cite{SuarezSerrato2026}.
Both sides of inequality \eqref{eq:rcd-entropy-degree-ineq} are invariant under rescaling the metric of $X$, and an $\RCD(K,N)$ space is an $\RCD(K',N)$ space for every $K' \leq K$.
Hence we may assume $K = -(N-1)$.

Let us include a brief synthesis of what is now called the barycenter method.
Douady and Earle \cite{Douady-Earle-1986} introduced the conformal barycenter of measures to extend circle homeomorphisms to the disc.
Besson, Courtois, and Gallot \cite{Besson-Courtois-Gallot:95} developed the barycenter method on smooth manifolds to prove the entropy rigidity of locally symmetric spaces of negative curvature.
Subsequent developments include the work of Bessi\`eres \cite{Bessieres-1998}, Connell and Farb \cite{Connell-Farb-2003}, Sambusetti \cite{Sambusetti-1999}, and Bessi\`eres, Besson, Courtois, and Gallot \cite{BBCG-2012}.
Storm \cite{Storm-2002, Storm-2006, Storm-2007} adapted the method to singular spaces, including Alexandrov spaces and cone-manifolds, in the setting of homotopy equivalences.
The version of the method used here relies on Sturm's $d^2$-barycenter on $\CAT(0)$ spaces \cite{Sturm03}.

In joint work with Connell, Dai, N\'u\~nez-Zimbr\'on, Perales, and Wei \cite{RCDbary}, we extended the barycenter method to non-collapsed $\RCD(K,N)$ spaces without boundary.
A degree theory for maps from $\RCD$ spaces was not available, so the lower bound in our previous work is expressed in terms of a homological index $\ind_H(f)$ and a fundamental index $\ind_\pi(f)$ \cite[Theorem 1.1]{RCDbary}.
These are defined through singular homology and fundamental groups, and when $X$ is a closed oriented manifold, $\ind_H(f) = |\deg(f)|$ \cite[Section 1]{RCDbary}.
The case of Alexandrov spaces with curvature bounded below was  treated after we showed that the analytical and topological degrees coincide \cite[Theorem 5.3]{SuarezSerrato2026}.
That proof uses the integral current structure of Alexandrov spaces constructed by Jaramillo, Perales, Rajan, Searle, and Siffert \cite{jaramillo2021alexandrov}, and Mitsuishi's isomorphism between integral current homology and singular homology \cite{mitsuishi2019homologies}, which applies to Alexandrov spaces by the local Lipschitz contractibility theorem of Mitsuishi and Yamaguchi \cite{mitsuishi2014llc}.
These two topological tools are not yet available for general non-collapsed $\RCD$ spaces.

The theory of currents in metric spaces was developed by Ambrosio and Kirchheim \cite{ambrosio2000currents} and, in a local form for locally compact spaces, by Lang \cite{Lang}.
Orientability in Ricci limit geometry was studied by Matveev and Portegies \cite{MatveevPortegies17} through intrinsic flat limits, and by Honda \cite{Honda17} through Sobolev differential forms.
For non-collapsed $\RCD$ spaces without boundary, Brena, Bru\`e, and Pigati \cite{BrenaBruePigati2024} proved the equivalence between topological orientability of the manifold part, Honda's formulation, and a version in terms of metric currents.
Mondino and Perales \cite{MP26} related orientability to local flat convergence.
Basso, Marti, and Wenger \cite{Basso-Marti-Wenger-2025} constructed metric fundamental classes for closed oriented metric manifolds under linear local contractibility.

Recently, Cavallucci, Marti, Mondino, and Perales \cite{CMMP2026} developed a theory of orientation by locally integral currents for metric manifolds and for singular spaces whose manifold part is dense and of full measure.
For a non-collapsed $\RCD(K,N)$ space $X$ with orientable manifold part, they constructed an orientation current $T$ with $\|T\| = \Hn$, unique up to sign, and identified its boundary with the orientation current of $\partial X$ \cite[Theorem 1.2]{CMMP2026}.
When $X$ is closed, $T$ is an integral cycle generating the top-dimensional integral current homology of $X$ \cite[Theorem 8.7]{CMMP2026}.

In this paper, we combine these two developments.
\Cref{prop:degree-equivalence} shows that for a Lipschitz map $f$ from a closed oriented non-collapsed $\RCD(K,N)$ space to a closed oriented $N$-manifold, the push-forward $f_\# T$ equals $\deg(f)$ times the fundamental class of the target.
That proof uses the representation of $T$ by local degrees on regular disks of $\mathcal{M}(X)$ \cite[Proposition 6.2 and Theorem 6.5]{CMMP2026}, together with the fact that the singular set $\mathcal{S}(X)$ has vanishing $\mathcal{H}^{N-1}$-measure when $\partial X = \emptyset$.
Inserting this identification into our previously developed barycenter method for $\RCD$ spaces \cite{RCDbary} proves \Cref{thm:rcd-entropy-degree}.

The inequalities of the barycenter method compare metric Jacobians defined with respect to the reference measure of the domain (cf \cite{SuarezSerrato2026}), and the argument requires this reference measure to be the Hausdorff measure.
For non-collapsed $\RCD(K,N)$ spaces, this holds by definition.

\medskip

While this manuscript was in its final stages of preparation, the preprint of Li, Pasqualetto, and Qian \cite{Li-Pasqualetto-Qian-2026} appeared. They successfully generalized the closed $\RCD$ space Lipschitz-volume rigidity theorem \cite{RCDbary} to non-collapsed $\RCD(K,N)$ spaces with boundary. For closed spaces, their result (\cite[Corollary 1.4]{Li-Pasqualetto-Qian-2026}) provides a concurrent, independent reinforcement of the metric rigidity engine used in the proof of the equality case of \Cref{thm:rcd-entropy-degree}.
Moreover, we remark that some technical simplifications to the arguments used here may become available if/once the recent preprint on bi-Lipschitz coordinates for non-collapsed $\RCD$ spaces \cite{OpenAI-2026} is independently verified by human readers, as is the standard for mathematical results to be accepted.

\subsection*{Organization of the paper}
\Cref{sec:orientation-degree} recalls the structure of non-collapsed $\RCD(K,N)$ spaces and the orientation currents of \cite{CMMP2026}, defines the degree, and proves its equivalence with the current-theoretic multiplicity.
\Cref{sec:barycenter} recalls the natural maps of the barycenter method on $\RCD$ spaces from \cite{RCDbary}, and the relation between metric Jacobians and the mass of push-forward currents.
\Cref{sec:proof} contains the proof of \Cref{thm:rcd-entropy-degree} and remarks on the orientability hypothesis, the Alexandrov case, and the equality case.


\section{Orientation currents and degree on non-collapsed \texorpdfstring{$\RCD$}{RCD} spaces}\label{sec:orientation-degree}

\subsection{Structure of non-collapsed \texorpdfstring{$\RCD$}{RCD} spaces}\label{subsec:ncRCD-structure}

Let $(X,\dist,\Hn)$ be a non-collapsed $\RCD(K,N)$ space.
Following standard terminology, we denote by $\mathcal{R}_N(X)$ the set of points at which every tangent cone is isometric to $\R^N$, and by $\mathcal{S}^k(X)$, $0 \leq k \leq N-1$, the set of points at which no tangent cone splits off a factor $\R^{k+1}$ \cite{MondinoNaber19, DPG}.
This stratification goes back to Cheeger and Colding \cite{Cheeger-Colding-1997} for Ricci limit spaces, and De Philippis and Gigli proved that $\dim_{\mathcal{H}} \mathcal{S}^k(X) \leq k$ \cite[Theorem 1.8]{DPG}.
The boundary of $X$ is $\partial X := \overline{\mathcal{S}^{N-1}(X) \setminus \mathcal{S}^{N-2}(X)}$ \cite{DPG}.
Bru\`e, Naber, and Semola proved that $\partial X$ is $(N-1)$-rectifiable and that the tangent cone at every point of $\mathcal{S}^{N-1}(X) \setminus \mathcal{S}^{N-2}(X)$ is unique and isometric to the half-space $\R^N_+$ \cite[Theorems 1.4 and 6.6]{BrueNaberSemola20}; see also \cite[Section 8]{CMMP2026}.
Kapovitch and Mondino \cite{KapMonGT} introduced a notion of boundary through iterated tangent cones and studied its topology.

Denote by $\mathcal{M}(X)$ the set of points of $X$ that have a neighbourhood homeomorphic to $\R^N$, and set $\mathcal{S}(X) := X \setminus \mathcal{M}(X)$.
A locally compact metric space is \emph{$N$-regular} if $\Hn$ is locally finite, $\mathcal{M}(X)$ is a degree-regular metric $N$-manifold in the sense of \cite[Definition 5.9]{CMMP2026}, $\mathcal{M}(X)$ is dense, and $\Hn(\mathcal{S}(X)) = 0$ \cite[Section 7]{CMMP2026}.
The following result synthesizes methods we will need in our arguments below.

\begin{lemma}\label{lem:ncRCD-structure}
Let $(X,\dist,\Hn)$ be a closed non-collapsed $\RCD(K,N)$ space, $N \geq 2$. Then:
\begin{enumerate}[label=(\roman*)]
    \item $X$ is $N$-rectifiable, and for $\Hn$-almost every $x \in X$ the approximate tangent space of $X$ at $x$ is isometric to Euclidean $\R^N$ \cite{MondinoNaber19, KellMondino, GigliPasqualetto:CAG}, \cite[Theorem 8.1]{CMMP2026};
   \item $X$ is $N$-regular, $\mathcal{R}_N(X) \subseteq \mathcal{M}(X)$, and $(\mathcal{M}(X),\dist)$ is a length space; in particular $\mathcal{M}(X)$ is connected \cite[Theorem 8.1]{CMMP2026};
    \item $\mathcal{S}(X) \subseteq \mathcal{S}^{N-2}(X)$. In particular, $\dim_{\mathcal{H}} \mathcal{S}(X) \leq N-2$ and 
    \[ \mathcal{H}^{N-1}(\mathcal{S}(X)) = 0.
    \]
\end{enumerate}
\end{lemma}

\begin{proof}
Items (i) and (ii) are contained in \cite[Theorem 8.1]{CMMP2026}, whose proof uses the rectifiability results of \cite{MondinoNaber19, KellMondino, GigliPasqualetto:CAG} and the regular coordinates of Mondino and Perales \cite[Proposition 3.8]{MP26}.

For (iii), by (ii) we have $\mathcal{S}(X) \subseteq X \setminus \mathcal{R}_N(X)$, and it suffices to show that $X \setminus \mathcal{R}_N(X) \subseteq \mathcal{S}^{N-2}(X)$.

We use two facts from the theory of non-collapsed spaces \cite[Theorems 1.6 and 1.10]{DPG}, see also \cite[Section 2.4]{BrueNaberSemola20}:

\begin{enumerate}
\item The density $\theta(x) := \lim_{r \to 0} \Hn(B(x,r))/(\omega_N r^N)$ exists at every point $x$ in $X$;
\item Every tangent cone $C$ at $x$ is a metric cone with vertex $o$ satisfying $\Hn(B(o,1)) = \omega_N \theta(x)$.
\end{enumerate}
 
Moreover, by the volume rigidity theorem \cite[Theorem 1.6 and Corollary 2.14]{DPG}, a tangent cone with $\Hn(B(o,1)) = \omega_N$ (equivalently, $\theta(x) = 1$) is isometric to $\R^N$.

Let $x$ be a point that is not in $\mathcal{S}^{N-1}(X)$.
Then some tangent cone at $x$ splits off $\R^N$, hence is isometric to $\R^N$, and $\theta(x) = 1$.
Consequently every tangent cone at $x$ is isometric to $\R^N$, that is, $x \in \mathcal{R}_N(X)$.

This shows $X \setminus \mathcal{R}_N(X) \subseteq \mathcal{S}^{N-1}(X)$.
Since $\partial X = \emptyset$, the set $\mathcal{S}^{N-1}(X) \setminus \mathcal{S}^{N-2}(X)$ is empty, so $\mathcal{S}^{N-1}(X) = \mathcal{S}^{N-2}(X)$.
The dimension bound follows from \cite[Theorem 1.8]{DPG}, and it implies $\mathcal{H}^{N-1}(\mathcal{S}(X)) = 0$.
\end{proof}

\subsection{Metric currents}\label{subsec:currents}

We use the theory of currents in metric spaces of Ambrosio and Kirchheim \cite{ambrosio2000currents}, in the local formulation of Lang \cite{Lang} adopted by Cavallucci--Marti--Mondino--Perales \cite{CMMP2026}.
Recall that a $k$-dimensional current $T$ on a locally compact metric space $Z$ is a multilinear functional on $(k+1)$-tuples $(g, \tau_1, \dots, \tau_k)$ of locally Lipschitz functions on $Z$, with $g$ compactly supported, satisfying continuity, locality, and the alternating property \cite[Definition 2.1]{Lang}, \cite[Section 3.1]{CMMP2026}.
A $k$-current $T$ is said to have \emph{locally finite mass} if there exists a locally finite Borel regular measure $\mu$ on $Z$ such that
\begin{equation}\label{eq:mass-definition}
|T(g, \tau_1, \dots, \tau_k)| \;\leq\; \prod_{j=1}^k \Lip(\tau_j|_{\spt(g)}) \int_Z |g|\, d\mu
\end{equation}
for every compactly supported Lipschitz function $g \colon Z \to \R$, and all $1$-Lipschitz functions $\tau_1, \dots, \tau_k$ on $\spt(g)$.
The minimal Borel regular measure $\mu$ satisfying \eqref{eq:mass-definition} is called the \emph{mass measure} of $T$ and is denoted by $\|T\|$ \cite[Theorem 4.3]{Lang}, \cite[Definition 3.10]{CMMP2026}.
The \emph{total mass} of $T$ is defined as $\bM(T) := \|T\|(Z) \in [0,\infty]$, and the \emph{support} of $T$ is $\spt(T) := \spt(\|T\|)$.

For a locally compact metric space $Z$, we denote by $\cI_{k,\mathrm{loc}}(Z)$ the locally integer rectifiable $k$-currents, and by $\bI_{k,\mathrm{loc}}(Z)$ the locally integral $k$-currents.
Namely those $T$ in $\cI_{k,\mathrm{loc}}(Z)$ whose boundary $\partial T$, defined by $\partial T(g, \tau_1, \dots, \tau_{k-1}) := T(1, g, \tau_1, \dots, \tau_{k-1})$, also has locally finite mass \cite[Sections 3 and 4]{CMMP2026}.
If $Z$ is compact, every current with locally finite mass has finite total mass $\bM(T) < \infty$. 
The classes of currents with finite mass, integer rectifiable currents, and integral currents in the sense of Lang coincide with those of Ambrosio and Kirchheim \cite[Sections 3.2 and 4.2]{CMMP2026}. 
We then write $\bI_k(Z)$.
A Lipschitz map $\varphi \colon Z \to W$ between compact metric spaces induces a push-forward $\varphi_\# \colon \bI_k(Z) \to \bI_k(W)$, given by $\varphi_\# T(g, \tau_1, \dots, \tau_k) := T(g \circ \varphi, \tau_1 \circ \varphi, \dots, \tau_k \circ \varphi)$, which commutes with the boundary operator and satisfies $\bM(\varphi_\# T) \leq \Lip(\varphi)^k \bM(T)$.

The integral current homology $H^{\mathrm{IC}}_k(Z)$ is the homology of the chain complex $(\bI_{\bullet,\mathrm{loc}}(Z),\partial)$ \cite[Section 4]{CMMP2026}; see also \cite{mitsuishi2019homologies, jaramillo2021alexandrov}.
If $\mathcal{H}^{k+1}(Z) = 0$, then $\bI_{k+1,\mathrm{loc}}(Z) = 0$ and $H^{\mathrm{IC}}_k(Z)$ is the group of $k$-cycles \cite[Lemma 4.12]{CMMP2026}.
This is based on the fact that the mass measure of a current $S$ with locally finite mass satisfies $\|S\| \ll \mathcal{H}^k$ \cite[Lemma 3.20]{CMMP2026}.

As a guiding example, in a closed oriented Riemannian $N$-manifold $Y$, the fundamental class $\cur{Y}$ in $\bI_N(Y)$ is given by 
\[\cur{Y}(g,\tau) = \int_Y g \det(D\tau)\, d\Hn. \]
It satisfies $\partial \cur{Y} = 0$ and $\|\cur{Y}\| = \Hn$ \cite[Section 4.2]{CMMP2026}.
The constancy theorem \cite[4.1.31]{Federer} (see also \cite[Corollary 6.8]{CMMP2026}) implies every $N$-cycle in $Y$ is an integer multiple of $\cur{Y}$.
In particular $\bM(k\cur{Y}) = |k|\, \Hn(Y)$.

\subsection{Orientation currents}\label{subsec:orientation-currents}

An $N$-regular space $X$ is \emph{orientable} if $\mathcal{M}(X)$ is orientable, and an \emph{orientation} of $X$ is a current $T \in \cI_{N,\mathrm{loc}}(X)$ with 
$ \|T\| = \lambda_E \Hn \mres E$ and $\spt(\partial T) \subseteq \mathcal{S}(X)$.
Here $E$ is the rectifiable part of $X$ and $\lambda_E$ its tangential function \cite[Definition 7.3 and Remark 7.4]{CMMP2026}.
Every orientation of $\mathcal{M}(X)$ extends uniquely to an orientation of $X$ \cite[Corollary 7.2]{CMMP2026}.
An $N$-regular space is \emph{locally uniquely orientable} if every point has arbitrarily small uniquely orientable ball neighbourhoods \cite[Definition 7.15]{CMMP2026}.

Let $D \subseteq \mathcal{M}(X)$ be a  a regular disk (a subset homeomorphic to a closed Euclidean ball with $\Hn(D) < \infty$ and $\Hn(\partial D) = 0$, \cite[Definition 5.1]{CMMP2026}) with an orientation $\mathfrak{o}$ of its interior. Cavallucci, Marti, Mondino, and Perales construct a current $T_D^{\mathfrak{o}}$ in $\cI_N(D)$ from local degrees of bi-Lipschitz charts, and prove that
\begin{equation}\label{eq:local-degree-formula}
    T_D^{\mathfrak{o}}(g,\pi) = \int_{\R^N} \sum_{x \in \pi^{-1}(z) \cap \mathrm{int}(D)} g(x)\, \deg(\pi,x)\, d\mathcal{H}^N(z)
\end{equation}
for every Lipschitz pair $(g,\pi)$ on $D$ such that $\Hn(\pi(P_D)) = 0$.
Here $P_D$ is the $N$-unrectifiable part of $D$ and $\deg(\pi,x)$ is the local degree of $\pi \colon \mathrm{int}(D) \to \R^N$ at $x$ \cite[Proposition 6.2]{CMMP2026}.

The following statement collects the results that we use.

\begin{theorem}[Cavallucci--Marti--Mondino--Perales \cite{CMMP2026}]\label{thm:CMMP}
Let $(X,\dist,\Hn)$ be a closed non-collapsed $\RCD(K,N)$ space, $N \geq 2$, with $\mathcal{M}(X)$ orientable.
\begin{enumerate}[label=(\roman*)]
    \item $X$ is locally uniquely orientable and uniquely orientable: it has exactly two orientations $\pm T$, and $\|T\| = \Hn$ \cite[Proposition 8.5]{CMMP2026}.
    \item $\partial T = 0$, so that $T$ lies in $\bI_N(X)$ \cite[Theorem 8.6]{CMMP2026}, and $H^{\mathrm{IC}}_N(X) = \Z\, T \cong \Z$ \cite[Theorem 8.7]{CMMP2026}.
    \item The restriction $T|_{\mathcal{M}(X)}$ determines a topological orientation of $\mathcal{M}(X)$ such that $T|_D = T_D^{\mathfrak{o}}$ for every regular disk $D \subseteq \mathcal{M}(X)$, where $\mathfrak{o}$ is the induced orientation of $\mathrm{int}(D)$ \cite[Theorem 6.5 and Remark 7.4]{CMMP2026}.
\end{enumerate}
\end{theorem}

By \Cref{lem:ncRCD-structure}(i) $X$ is $N$-rectifiable, so the unrectifiable set of every regular disk is $\Hn$-negligible.
Therefore, Equation \eqref{eq:local-degree-formula} holds for every Lipschitz pair $(g,\pi)$ on $D$.
In the rest of the paper, an \emph{oriented} closed non-collapsed $\RCD(K,N)$ space is one with orientable $\mathcal{M}(X)$ together with a choice of orientation $T$, and $\mathcal{M}(X)$ carries the topological orientation given by \Cref{thm:CMMP}(iii).

\subsection{Degree}\label{subsec:degree}

Throughout this subsection, let $(X,\dist,\Hn)$ be an oriented, closed, and non-collapsed $\RCD(K,N)$ space with orientation $T$, $N \geq 2$, and $Y$ a closed oriented Riemannian $N$-manifold.
We use the local degree and the degree of proper maps between oriented topological manifolds, following Dold \cite[Chapter VIII]{Dold} (ap. \cite[Section 5.1]{CMMP2026}).

\begin{lemma}\label{lem:regular-values}
Let $f \colon X \to Y$ be Lipschitz, and set $\Sigma_f := f(\mathcal{S}(X))$, $V_f := Y \setminus \Sigma_f$, and $U_f := f^{-1}(V_f)$.
Then $\Sigma_f$ is compact with $\mathcal{H}^{N-1}(\Sigma_f) = 0$, the set $V_f$ is open, dense, and connected, $U_f$ is an open subset of $\mathcal{M}(X)$, and $f|_{U_f} \colon U_f \to V_f$ is proper.
\end{lemma}

\begin{proof}
The set $\mathcal{S}(X)$ is closed in the compact space $X$, so $\Sigma_f$ is compact.
By \Cref{lem:ncRCD-structure}(iii), $\mathcal{H}^{N-1}(\Sigma_f) \leq \Lip(f)^{N-1} \mathcal{H}^{N-1}(\mathcal{S}(X)) = 0$.
A theorem of Szpilrajn implies that a metric space $A$ with $\mathcal{H}^{N-1}(A) = 0$ has topological dimension at most $N-2$ \cite[Chapter VII, Theorem 2]{Hurewicz-Wallman}.
Recall that a closed subset of $\R^N$ of topological dimension at most $N-2$ has empty interior and does not separate any connected open subset of $\R^N$ \cite[Chapter IV, Theorem 4]{Hurewicz-Wallman}.
Covering the connected manifold $Y$ by coordinate balls, se see that the same holds for closed subsets of $Y$.
Hence $V_f$ is open, dense, and connected.
By construction $U_f \cap \mathcal{S}(X) = \emptyset$, so $U_f \subseteq \mathcal{M}(X)$, and $U_f$ is open because $f$ is continuous.
If $C \subseteq V_f$ is compact, then $f^{-1}(C)$ is a closed subset of the compact space $X$ contained in $U_f$, hence compact.
\end{proof}

As $U_f$ is an open subset of the oriented manifold $\mathcal{M}(X)$, and $f|_{U_f}$ is proper, the degree $\deg_y(f|_{U_f})$ (in $\Z$) over a point $y$ in $V_f$ is well defined.
Moreover, it is locally constant in $y$ \cite[Chapter VIII, Section 4]{Dold}.
The set $V_f$ is connected, so it does not depend on the point $y$ in $V_f$.

\begin{definition}[Degree of a Lipschitz map]\label{def:deg-Lipschitz}
For a Lipschitz map $f \colon X \to Y$, the \emph{degree} of $f$ is
\[
\deg(f) := \deg_y\big(f|_{U_f}\big), \qquad y \in V_f .
\]
\end{definition}

If $X$ is a closed oriented topological manifold, then $\mathcal{S}(X) = \emptyset$, $U_f = X$, $V_f = Y$, and \Cref{def:deg-Lipschitz} is the classical degree, characterized by $f_*[X] = \deg(f)[Y]$ in singular homology \cite[Chapter VIII, Section 4]{Dold}.
By the additivity of the local degree \cite[Chapter VIII, Proposition 4.7]{Dold}, for every $y$ in $V_f$ with finite $f^{-1}(y)$,
\begin{equation}\label{eq:degree-additivity}
\deg(f) = \sum_{x \in f^{-1}(y)} \deg(f,x).
\end{equation}

Recall that by the general Eilenberg--coarea inequality \cite[2.10.25]{Federer}, if $\varphi \colon Z \to W$ is a Lipschitz map between metric spaces and $m \geq k \geq 0$, then for every Borel set $A \subseteq Z$,
\begin{equation}\label{eq:coarea-inequality}
\int_W^* \mathcal{H}^{m-k}\big(A \cap \varphi^{-1}(w)\big)\, d\mathcal{H}^k(w) \;\leq\; \frac{\omega_{m-k}\,\omega_k}{\omega_m}\, \Lip(\varphi)^k\, \mathcal{H}^m(A).
\end{equation}
Here, $\int_W^*$ denotes the upper integral and $\omega_d$ is the Lebesgue measure of the unit ball in $\R^d$ (with $\omega_0 = 1$).
This estimate specializes in the equal-dimension case $m = k = N$, where $\mathcal{H}^0 = \#$ is the counting measure and the dimensional constant equals $1$.
Then \eqref{eq:coarea-inequality} takes the following form (\cite[Lemma 2.1]{CMMP2026}),
\[
\int_W^* \#\big(A \cap \varphi^{-1}(w)\big)\, d\mathcal{H}^N(w) \;\leq\; \Lip(\varphi)^N\, \mathcal{H}^N(A).
\]
This implies that if $\mathcal{H}^N(A) < \infty$, the pre-image $A \cap \varphi^{-1}(w)$ is finite for $\mathcal{H}^N$-almost every $w$ in $W$ (and $\mathcal{H}^N(\varphi(A)) = 0$ whenever $\mathcal{H}^N(A) = 0$).

The next result is the key linking the local analytic degree in \Cref{def:deg-Lipschitz} with the fundamental current.

\begin{proposition}[Analytic and topological degree]\label{prop:degree-equivalence}
Let $f \colon X \to Y$ be a Lipschitz map. Then
\[
f_\# T = \deg(f)\, \cur{Y}.
\]
In particular, 
\begin{equation}\label{eq:push-deg-curr}
 \deg(f) = \frac{1}{\Hn(Y)} \int_Y \theta_{f_\# T}\, d\Hn, 
\end{equation}
where $\theta_{f_\#T}$ is the signed multiplicity of $f_\#T = \theta_{f_\# T}\cur{Y}$ with respect to the orientation of $Y$.
\end{proposition}

\begin{proof}
By \Cref{thm:CMMP}(ii), $T$ lies in $\bI_N(X)$ and $\partial T = 0$.
Hence $f_\# T \in \bI_N(Y)$ and 
\[
\partial(f_\# T) = f_\#(\partial T) = 0.
\]
The constancy theorem on $Y$ \cite[4.1.31]{Federer} (see also \cite[Corollary 6.8]{CMMP2026} and \Cref{subsec:currents}) implies that there exists $k \in \Z$ such that $f_\# T = k \cur{Y}$.
It remains to show that $k = \deg(f)$.

Let $V_f$ and $U_f = f^{-1}(V_f)$ be as in \Cref{lem:regular-values}, and fix a point $y_0 \in V_f$.
Consider a closed geodesic ball $B \subseteq V_f$ centered at $y_0$ with radius smaller than the injectivity radius of $Y$.
Let $\phi \colon B \to \R^N$ be a smooth orientation-preserving normal coordinate chart onto a closed Euclidean ball (which is, in particular, a bi-Lipschitz homeomorphism), and denote its inverse by $\psi \colon \phi(B) \to B$.
By McShane's extension theorem \cite[Theorem 1]{McShane1934}, we may extend $\phi$ to a Lipschitz map $\phi \colon Y \to \R^N$ (noting that the extension need not be injective outside $B$).
Let $g \colon Y \to [0,\infty)$ be a Lipschitz function, not identically zero, with compact support $\spt(g) \subseteq \mathrm{int}(B)$. 
Since $\phi|_B$ is orientation-preserving, $\det(D\phi) = |\det(D\phi)| > 0$ almost everywhere on $B$. 

On the one hand, the area formula applied to $\phi|_B$ yields
\begin{equation}\label{eq:Y-side}
k\cur{Y}(g,\phi) = k\int_B g(y) \det(D\phi(y))\, d\Hn(y) = k\int_{\phi(B)} g(\psi(z))\, d\mathcal{H}^N(z),
\end{equation}
where the last integral is strictly positive by the choice of $g$.
On the other hand, by the definition of the push-forward current, 
\[ f_\# T(g,\phi) = T(g\circ f, \pi),\]
where $\pi := \phi \circ f \colon X \to \R^N$.
By \Cref{lem:regular-values}, $f|_{U_f} \colon U_f \to V_f$ is proper, so $C := f^{-1}(\spt(g))$ is a compact subset of $U_f \subseteq \mathcal{M}(X)$.
Because $\mathcal{M}(X)$ is a degree-regular metric manifold (\Cref{lem:ncRCD-structure}(ii)), ther exist finitely many regular disks $D_1, \dots, D_m \subseteq f^{-1}(\mathrm{int}(B)) \subseteq U_f$ whose interiors cover $C$ \cite[Lemma 5.2]{CMMP2026}.
Using standard Lipschitz cutoff functions on $X$, we may choose Lipschitz functions $\chi_1,\dots,\chi_m \colon X \to [0,1]$ such that $\spt(\chi_j) \subseteq \mathrm{int}(D_j)$ for each $j$ and $\sum_{j=1}^m \chi_j = 1$ on $C$.

As $g \circ f$ vanishes outside $C$, we have $g \circ f = \sum_{j=1}^m \chi_j \cdot (g \circ f)$ on $X$.
By the linearity and locality of $T$, together with \Cref{thm:CMMP}(iii) and the local degree formula \eqref{eq:local-degree-formula}, we obtain
\begin{align*}
T(g\circ f,\pi) &= \sum_{j=1}^m T_{D_j}^{\mathfrak{o}}\big(\chi_j \cdot (g\circ f), \pi|_{D_j}\big) \\
&= \sum_{j=1}^m \int_{\R^N} \sum_{x \in \pi^{-1}(z) \cap \mathrm{int}(D_j)} \chi_j(x)\, g(f(x))\, \deg(\pi,x)\, d\mathcal{H}^N(z) \\
&= \int_{\R^N} \sum_{x \in \pi^{-1}(z) \cap C} g(f(x))\, \deg(\pi,x)\, d\mathcal{H}^N(z).
\end{align*}
The last equality uses that $\spt(g \circ f) \subseteq C$ and $\sum_{j=1}^m \chi_j = 1$ on $C$.

By the coarea inequality \eqref{eq:coarea-inequality} applied with $m = k = N$ (see \cite[2.10.25]{Federer} and \cite[Lemma 2.1]{CMMP2026}), since $\Hn(X) < \infty$, the fiber $f^{-1}(y)$ is finite for $\Hn$-almost every $y \in B$, and thus $\pi^{-1}(z) \cap f^{-1}(B)$ is finite for $\mathcal{H}^N$-almost every $z \in \phi(B)$.
Notice that if $x \in C$, then $f(x) \in \spt(g) \subseteq \mathrm{int}(B)$, so $\pi(x) = \phi(f(x)) \in \phi(\spt(g)) \subseteq \phi(B)$.
Conversely, because $\phi$ is injective on $B$, for any $z \in \phi(B)$ and any $x \in f^{-1}(B)$, the condition $\pi(x) = z$ is equivalent to $f(x) = \psi(z)$.
The vanishing $g(f(x)) = 0$ for $x \in f^{-1}(\psi(z)) \setminus C$ implies we may replace the summation domain $\pi^{-1}(z) \cap C$ by the pre-image $f^{-1}(\psi(z))$ for every point $z$ in $\phi(B)$.

Observe that $\deg(\phi, f(x)) = 1$, because for every $x \in f^{-1}(\mathrm{int}(B)) \subseteq U_f$, the map $\phi$ is an orientation-preserving homeomorphism on a neighbourhood of $f(x)$.
The multiplicativity of the local degree \cite[Chapter VIII, Corollary 4.6]{Dold} thus implies
\[
\deg(\pi,x) = \deg(\phi,f(x))\,\deg(f,x) = \deg(f,x).
\]
Substituting this into our integral and applying the additivity formula \eqref{eq:degree-additivity} at $\mathcal{H}^N$-almost every $z \in \phi(B)$ (with $y = \psi(z) \in V_f$), we find
\begin{align*}
f_\# T(g,\phi) &= \int_{\phi(B)} g(\psi(z)) \sum_{x \in f^{-1}(\psi(z))} \deg(f,x)\, d\mathcal{H}^N(z) \\
&= \deg(f) \int_{\phi(B)} g(\psi(z))\, d\mathcal{H}^N(z).
\end{align*}
Comparing with \Cref{eq:Y-side}, and dividing by the positive integral  factor on the right hand side, we find  $k = \deg(f)$.

Finally, as $f_\# T = \deg(f)\cur{Y}$ and $\|\cur{Y}\| = \Hn$, the multiplicity of $f_\# T$ with respect to the orientation of $Y$ is constant and equal to $\theta_{f_\# T} = \deg(f)$ $\Hn$-almost everywhere.
From which we obtain the integral formula in \Cref{eq:push-deg-curr}.
\end{proof}
 
For completeness we next review how to extend the degree to continuous maps.

\begin{lemma}[Homotopy invariance]\label{lem:homotopy-invariance}
Every continuous map $f \colon X \to Y$ is homotopic to a Lipschitz map.
If $f_0, f_1 \colon X \to Y$ are Lipschitz and homotopic, then $f_{0\#} T = f_{1\#} T$; in particular $\deg(f_0) = \deg(f_1)$.
\end{lemma}

\begin{proof}
Let $Y \subseteq \R^m$ be an isometric embedding, $\mathcal{N}_\delta$ a tubular neighbourhood of $Y$ of radius $\delta>0$, and $r \colon \mathcal{N}_\delta \to Y$ the nearest-point retraction, which is Lipschitz.

Recall that Lipschitz functions are uniformly dense in the continuous functions on a compact metric space, by the Stone--Weierstrass theorem (see also \cite[Theorem 6.8]{Heinonen} or \cite[Section 2.7]{SuarezSerrato2026}).
So after composing uniform Lipschitz approximations with $r$ a continuous map from a compact metric space to $Y$ may be uniformly approximated by Lipschitz maps.
Namely, a pair of maps $f_0,f_1$ with $\sup_x |f_0(x) - f_1(x)| < \delta$ are joined by an homotopy $(t,x) \mapsto r((1-t)f_0(x) + tf_1(x))$, which is Lipschitz whenever $f_0$ and $f_1$ are.

Let $H \colon [0,1] \times X \to Y$ be a homotopy from $f_0$ to $f_1$.
Applying the preceding paragraph to $H$ on the compact space $[0,1] \times X$ gives a Lipschitz map $H'$ uniformly close to $H$, and concatenating $H'$ with Lipschitz homotopies from $f_0$ to $H'_0$ and from $H'_1$ to $f_1$ gives a Lipschitz homotopy $G \colon [0,1] \times X \to Y$ from $f_0$ to $f_1$.
Then, by the homotopy formula for metric currents \cite[Theorem 3.2]{Wenger-2005} (see also \cite[Eq. (4.4)]{Lang}),
\[
f_{1\#} T - f_{0\#} T = \partial\, G_\#\big(\cur{0,1} \times T\big) + G_\#\big(\cur{0,1} \times \partial T\big) .
\]
The second term vanishes because $\partial T = 0$.
The current $G_\#(\cur{0,1}\times T)$ has dimension $N+1$ and finite mass in $Y$, so its mass measure is absolutely continuous with respect to $\mathcal{H}^{N+1}$ \cite[Lemma 3.20]{CMMP2026} (see also \cite[Theorem 3.7]{ambrosio2000currents}).
As $\mathcal{H}^{N+1}(Y) = 0$, this current vanishes on $Y$.
Hence $f_{1\#}T = f_{0\#}T$, and by \Cref{prop:degree-equivalence} $\deg(f_0) = \deg(f_1)$.
\end{proof}

\begin{definition}[Degree of a continuous map]\label{def:deg-continuous}
The \emph{degree} $\deg(f)$ of a continuous map $f \colon X \to Y$ is the degree of any Lipschitz map homotopic to $f$.
\end{definition}

By \Cref{lem:homotopy-invariance}, \Cref{def:deg-continuous} is independent of the Lipschitz representative, it agrees with \Cref{def:deg-Lipschitz} for Lipschitz maps, and it is a homotopy invariant.
Replacing $T$ by $-T$, or the orientation of $Y$ by the opposite one, changes the sign of $\deg(f)$; the absolute value $|\deg(f)|$ does not depend on these choices.

\begin{corollary}\label{cor:surjective}
If $f \colon X \to Y$ is continuous and $\deg(f) \neq 0$, then $f$ is surjective.
\end{corollary}

\begin{proof}
If $f$ is Lipschitz, then by \Cref{prop:degree-equivalence} and the fact that $\spt\cur{Y} = Y$, we have
\[
Y = \spt\big(\deg(f)\cur{Y}\big) = \spt(f_\# T) \subseteq f(\spt T) \subseteq f(X).
\]
If $f$ is only continuous and there exists a point $y$ in $Y \setminus f(X)$, then because $f(X)$ is compact, there exists $\epsilon > 0$ such that the ball $B(y, 2\epsilon)$ is disjoint from $f(X)$.
By the proof of \Cref{lem:homotopy-invariance}, there exists a Lipschitz map $f' \colon X \to Y$ homotopic to $f$ satisfying $\sup_{x \in X} \dist(f(x), f'(x)) < \epsilon$.
By the triangle inequality, $y \notin f'(X)$, so $f'$ is not surjective.
However, by \Cref{def:deg-continuous}, $\deg(f') = \deg(f) \neq 0$, which contradicts the Lipschitz case.
\end{proof}


\section{The barycenter method on non-collapsed \texorpdfstring{$\RCD$}{RCD} spaces}\label{sec:barycenter}

\subsection{Natural maps}\label{subsec:natural-maps}

A construction of the natural maps of the barycenter method for $\RCD$ spaces, continuing from the work of  Besson--Courtois--Gallot \cite{Besson-Courtois-Gallot:95}, Sambusetti \cite{Sambusetti-1999}, and Sturm \cite{Sturm03},  is available\cite[Section 4]{RCDbary} (see also \cite[Section 6]{SuarezSerrato2026}).

Let $X$, $Y$, and $f$ be as in \Cref{thm:rcd-entropy-degree}, and let $\overline{X} = \til{X}/\ker f_*$ and $\Gamma = \pi_1(X)/\ker f_*$ be as in the introduction, with lift $\overline{f} \colon \overline{X} \to \til{Y}$ and covering map $p \colon \overline{X} \to X$.
We equip $\overline{X}$ with the length metric for which $p$ is a local isometry, and with the $\Gamma$-invariant lift $\overline{\mathfrak{m}}$ of $\Hn$, which is the $N$-dimensional Hausdorff measure of $\overline{X}$.
With this structure, $\overline{X}$ is a non-collapsed $\RCD(K,N)$ space \cite{Mondino-Wei-2019}, cf.\ \cite[Section 4]{RCDbary}.
When $K < 0$ is satisfied, the Bishop--Gromov inequality implies (cf \cite{RCDbary}), \[ h(\overline{X}) \leq \sqrt{(N-1)|K|} < \infty.\]

Let $s > h(\overline{X})$, and $x$ be a point  in $\overline{X}$.
Consider the finite measure  $\mu_x^s$ on $\overline{X}$ defined implicitly through its Radon--Nikodym derivative as follows,
\[
\frac{d\mu_x^s}{d\overline{\mathfrak{m}}}(z) = e^{-s\, \dist(x,z)}.
\]
The space $\til{Y}$ is a complete $\CAT(0)$ space, and every probability measure $\nu$ on $\til{Y}$ with finite first moment has a unique $d^2$-barycenter $\bary(\nu)$.
It is the minimum point of $z \mapsto \int (\dist(y,z)^2 - \dist(o,y)^2)\, d\nu(y)$ for any $o$ in $\til{Y}$ \cite[Proposition 4.3]{Sturm03}.
The natural map $\overline{F}_s \colon \overline{X} \to \til{Y}$ is defined as,
\[
\overline{F}_s(x) := \bary\left( \frac{\overline{f}_* \mu_x^s}{\mu_x^s(\overline{X})} \right).
\]

\begin{lemma}[{\cite[Lemmas 4.3, 4.5, and 4.7]{RCDbary}}]\label{lem:natural-maps}
For every $s > h(\overline{X})$, the map $\overline{F}_s$ is equivariant with respect to the homomorphism $\Gamma \to \pi_1(Y)$ induced by $f$, and descends to a Lipschitz map $F_s \colon X \to Y$ which is homotopic to $f$.
\end{lemma}

Moreover, by \Cref{lem:natural-maps} and \Cref{def:deg-continuous}, $\deg(F_s) = \deg(f)$ for every $s > h(\overline{X})$, and, by \Cref{prop:degree-equivalence},
\begin{equation}\label{eq:Fs-pushforward}
F_{s\#} T = \deg(f)\, \cur{Y}.
\end{equation}

\subsection{Jacobians and the mass of push-forwards}\label{subsec:jacobians}

Let $F \colon X \to Y$ be a Lipschitz map.
Recall from \Cref{lem:ncRCD-structure}(i) that $X$ is countably $\Hn$-rectifiable, meaning it can be covered up to an $\Hn$-negligible set by a countable family of bi-Lipschitz images $\varphi_i(E_i)$ of Borel subsets $E_i \subseteq \R^N$ \cite{MondinoNaber19, KellMondino}.
The differentiation theory of Kirchheim \cite[Corollary 5]{kirchheim1994} and Ambrosio--Kirchheim \cite[Theorem 8.1]{AK-rect} allows us to define the \emph{approximate tangent space} $\Tan(X,x)$ at $\Hn$-almost every $x$ in $X$. 
Specifically, at every density point $x = \varphi_i(z)$ of $E_i$ where the metric differential $\mathrm{md}_z \varphi_i$ exists and is non-degenerate, $\Tan(X,x)$ is given by the $N$-dimensional normed space $(\R^N, \mathrm{md}_z \varphi_i)$.
This space is independent of the choice of chart up to linear isometry.
For a non-collapsed $\RCD(K,N)$ space $X$, this approximate tangent space is isometric to Euclidean $\R^N$ at $\Hn$-almost every $x$ in $X$, where it coincides with the unique metric tangent cone $\lim_{r \to 0} (X, r^{-1}\dist, x)$ in the pointed Gromov--Hausdorff sense \cite{DPG, GigliPasqualetto:CAG}.

In addition, by Rademacher's theorem for rectifiable metric spaces \cite[Theorem 8.1]{AK-rect}, at $\Hn$-almost every $x$ in $X$ the map $F$ is differentiable, so it admits an approximate \emph{tangential differential}
\[
d_x F \colon \Tan(X,x) \longrightarrow T_{F(x)}Y.
\]
It is a linear map between $N$-dimensional Euclidean spaces satisfying
\[
\lim_{y \to x,\, y \in X'} \frac{\dist_Y\big(F(y), \exp_{F(x)}(d_x F(v_y))\big)}{\dist(x,y)} = 0
\]
along an $\Hn$-density-$1$ subset $X' \subseteq X$ at $x$ (where $v_y \in \Tan(X,x)$ is the local chart displacement corresponding to $y$).

A definition of the Jacobian $\Jac_x F$, in terms of the coarea factor of the metric differential, following Ambrosio--Kirchheim \cite{AK-rect} and Reichel \cite{Reichel-2009}, is also available \cite[Section 2.2]{RCDbary}.

\begin{lemma}[{\cite[Lemma 2.8]{RCDbary}}]\label{lem:Jacobian}
For $\Hn$-almost every $x \in X$,
\[
\Jac_x F = |\det d_x F| = \limsup_{r \to 0} \frac{\Hn(F(B(x,r)))}{\Hn(B(x,r))}.
\]
\end{lemma}

In the setting of Alexandrov spaces, Lytchak's differentiation theorem \cite{lytchak2005differentiation} was used  to identify the Jacobian of the area formula with the volume-ratio Jacobian \cite[Section 4]{SuarezSerrato2026}.
For non-collapsed $\RCD$ spaces, \Cref{lem:Jacobian} plays the same role, and it uses that the reference measure is $\Hn$.

\begin{lemma}[Mass of push-forwards]\label{lem:mass-pushforward}
For every Lipschitz map $F \colon X \to Y$,
\[
\bM(F_\# T) \leq \int_Y \#F^{-1}(y)\, d\Hn(y) = \int_X \Jac_x F\, d\Hn(x).
\]
\end{lemma}

\begin{proof}
The equality is the area formula for Lipschitz maps on rectifiable metric spaces \cite[Theorem 8.2]{AK-rect}, see also \cite[Theorem 7]{kirchheim1994}, combined with \Cref{lem:Jacobian}.
For the inequality, let $\{\rho_i, K_i, \theta_i\}$ be a parametrization of $T$ in the sense of \cite[Definition 4.6]{CMMP2026}, with $K_i \subseteq X$ compact and pairwise disjoint, $\rho_i \colon K_i \to \R^N$ bi-Lipschitz, $T = \sum_i (\rho_i^{-1})_\# \cur{\theta_i}$, and 
\[ \|T\| = \sum_i \|(\rho_i^{-1})_\# \cur{\theta_i}\| . \]

\Cref{thm:CMMP}(i) implies $\|T\| = \Hn$, and the tangential function of $X$ equals $1$ almost everywhere.
Hence by \Cref{lem:ncRCD-structure}(i), the mass representation \cite[Lemma 4.10]{CMMP2026} yields $|\theta_i| = 1$ almost everywhere on $\rho_i(K_i)$.
Set $\psi_i := F \circ \rho_i^{-1}$. Then $F_\# T = \sum_i \psi_{i\#} \cur{\theta_i}$, and the area formula for Lipschitz maps from subsets of $\R^N$ to $Y$ gives
\[
\bM\big(\psi_{i\#}\cur{\theta_i}\big) \leq \int_{\rho_i(K_i)} J\psi_i(z)\, d\mathcal{H}^N(z) = \int_{K_i} \Jac_x F\, d\Hn(x).
\]
Here, the equality follows from the chain rule for metric differentials and the area formula for $\rho_i^{-1}$ \cite{kirchheim1994, AK-rect}.
Summing over $i$ we obtain the inequality, as the sets $K_i$ cover $X$ up to an $\Hn$-null set.
\end{proof}

As usual when applying the barycenter method, the following estimate is key.

\begin{proposition}[{Entropy Jacobian bound, \cite[Proposition 4.8]{RCDbary}}]\label{prop:entropy-jacobian}
Let $N \geq 3$, and let $X$, $Y$, $f$ be as in \Cref{thm:rcd-entropy-degree}. For every $s > h(\overline{X})$ and $\Hn$-almost every $x \in X$,
\[
\Jac_x F_s \leq \left( \frac{s}{h(\til{Y})} \right)^N .
\]
\end{proposition}

\section{Proof of  \texorpdfstring{\Cref{thm:rcd-entropy-degree}}{the main theorem}}\label{sec:proof}

\begin{proof}[Proof of the inequality in \Cref{thm:rcd-entropy-degree}]
Fix an orientation $T$ of $X$ and an orientation of $Y$.
Let $s > h(\overline{X})$, and $F_s \colon X \to Y$ be the natural map of \Cref{lem:natural-maps}, which is Lipschitz and homotopic to $f$.
By Equation \eqref{eq:Fs-pushforward}, \Cref{lem:mass-pushforward}, and \Cref{prop:entropy-jacobian},
\begin{equation}\label{eq:main-chain}
\begin{aligned}
|\deg(f)|\, \Hn(Y) &= \bM\big(\deg(f)\cur{Y}\big) = \bM(F_{s\#}T) \\
&\leq \int_X \Jac_x F_s\, d\Hn(x) \leq \left( \frac{s}{h(\til{Y})} \right)^N \Hn(X).
\end{aligned}
\end{equation}
Letting $s \to h(\overline{X})$ we find inequality \eqref{eq:rcd-entropy-degree-ineq}.
\end{proof}

\begin{remark}[Orientability]\label{rmk:orientability}
If $\mathcal{M}(X)$ is not orientable, then $H^{\mathrm{IC}}_N(X) = \{0\}$ \cite[Theorem 8.7]{CMMP2026}, so $X$ does not admit a non-zero integral $N$-cycle, and the current-theoretic multiplicity of a map vanishes.
In contrast, the inequality obtained in our previous work on the barycenter method for $\RCD$ spaces  \cite[Theorem 1.1]{RCDbary}, stated with the indices $\ind_H(f)$ and $\ind_\pi(f)$, does not require orientability.
\end{remark}

\begin{remark}[Comparison with the homological index]\label{rmk:index-comparison}
When $X$ is a closed oriented topological manifold, \Cref{def:deg-continuous} coincides with the classical degree, and then $\ind_H(f) = |\deg(f)|$ \cite[Section 1]{RCDbary}.
More generally, suppose that Mitsuishi's homomorphism $\Psi \colon H^{\mathrm{IC}}_N(X) \to H_N(X;\Z)$  \cite{mitsuishi2019homologies} is defined (for instance, when $X$ is locally Lipschitz contractible).
Then $f_*\Psi([T]) = \deg(f)[Y]$ in $H_N(Y;\Z) \cong \Z[Y]$, so $\deg(f)\Z[Y] \subseteq f_*H_N(X;\Z)$.
By \cite[Equation (1.2)]{RCDbary}, the homological index is given by $\ind_H(f) = [H_N(Y;\Z) : f_*H_N(X;\Z)]$ when this subgroup index is finite, and $0$ otherwise.
Consequently, whenever $\deg(f) \neq 0$, the index $[H_N(Y;\Z) : f_*H_N(X;\Z)]$ is finite and $\ind_H(f)$ divides $|\deg(f)|$.
\end{remark}

\begin{remark}[Alexandrov spaces]\label{rmk:Alexandrov}
Let $X$ be a closed (compact, without boundary) $N$-dimensional Alexandrov space with curvature bounded below by $\kappa$.
Then $(X,\dist,\Hn)$ is a non-collapsed $\RCD((N-1)\kappa,N)$ space \cite{petrunin2011, zhang2010ricci}.
By Burago, Gromov, and Perelman \cite{burago1992ad}, the tangent cone at $x \in X$ is unique and equal to the Euclidean cone over the space of directions $\Sigma_x$, which has no boundary.
Hence its boundary (in the sense of De Phillipis--Gigli \cite{DPG}) is empty
If it splits off $\R^{N-1}$, then it is $\R^{N-1} \times C(\Sigma')$ with $\Sigma'$ a $0$-dimensional space of directions without boundary, that is, two points, and the tangent cone is $\R^N$.
Hence $\mathcal{S}^{N-1}(X) \setminus \mathcal{S}^{N-2}(X) = \emptyset$.
Orientability of $X$ in the sense of Mitsuishi \cite{mitsuishi2016orientability} is equivalent to orientability of $\mathcal{M}(X)$.
Indeed, by Mitsuishi's isomorphism \cite{mitsuishi2019homologies}, which applies because of local Lipschitz contractibility \cite{mitsuishi2014llc}, $H_N(X;\Z) \cong H^{\mathrm{IC}}_N(X)$.
If $\mathcal{M}(X)$ is orientable, then $H^{\mathrm{IC}}_N(X) \cong \Z$ by \Cref{thm:CMMP}(ii), and $X$ is orientable by \cite[Theorem 3.8]{SuarezSerrato2026}.
If $\mathcal{M}(X)$ is not orientable, then $H^{\mathrm{IC}}_N(X) = \{0\}$ by \cite[Theorem 8.7]{CMMP2026}, and $X$ is not orientable.
In the orientable case, the orientation current $T$ and the current $T_X$ of \cite{jaramillo2021alexandrov} both span $H^{\mathrm{IC}}_N(X)$.
Therefore $T = \pm T_X$, and \Cref{prop:degree-equivalence} together with \cite[Theorem 5.3]{SuarezSerrato2026} shows that \Cref{def:deg-continuous} agrees with the topological degree used in the case of Alexandrov spaces\cite{SuarezSerrato2026}.
This shows that \Cref{thm:rcd-entropy-degree} contains the corresponding inequality for Alexandrov spaces \cite[Theorem 1.1]{SuarezSerrato2026}.
\end{remark}

We now establish the rigidity in the equality case.

\begin{proof}[Proof of the equality case of \Cref{thm:rcd-entropy-degree}]
Assume that equality holds in inequality \eqref{eq:rcd-entropy-degree-ineq}, and set $k := |\deg(f)| = [\pi_1(Y) : f_*\pi_1(X)] \geq 1$. 

Let $p \colon \hat{Y} \to Y$ be the finite Riemannian covering space associated with the subgroup $f_*\pi_1(X) \leq \pi_1(Y)$. 
By classical covering space theory, $p$ is a local isometry of degree $\deg(p) = k$.
It follows that $\Hn(\hat{Y}) = k \Hn(Y)$, and the universal cover of $\hat{Y}$ is isometric to $\til{Y}$, implying $h(\tilde{\hat{Y}}) = h(\til{Y})$.
The map $f$ lifts to a continuous map $\hat{f} \colon X \to \hat{Y}$ such that $f = p \circ \hat{f}$. 
By the multiplicativity of the degree (\cite[Chapter VIII, Corollary 4.6]{Dold}), we have $|\deg(f)| = |\deg(p)| \cdot |\deg(\hat{f})|$, which implies $|\deg(\hat{f})| = 1$. 
Observe that the intermediate cover $\overline{X}$ associated to $\hat{f}$ is exactly the same as the one associated to $f$, because $f_*\pi_1(X) = p_*\pi_1(\hat{Y})$.

The equality hypothesis on inequality \eqref{eq:rcd-entropy-degree-ineq} can thus be rewritten as 
\[h(\overline{X})^N \Hn(X) = h(\tilde{\hat{Y}})^N \Hn(\hat{Y}).\]
Recall that under a homothetic rescaling of the metric by a factor $c > 0$, the volume $\Hn(X)$ scales by $c^N$, while the entropy $h(\overline{X})$ scales by $c^{-1}$.
Consequently, the product $h(\overline{X})^N \Hn(X)$ is a scale-invariant quantity.
By suitably rescaling $X$, we may then assume (without loss of generality) that $h(\overline{X}) = h(\tilde{\hat{Y}}) = h(\til{Y})$.
The equality hypothesis then simplifies to $\Hn(X) = \Hn(\hat{Y})$.

We apply the barycenter method directly to the degree-one map $\hat{f} \colon X \to \hat{Y}$. 
Following the construction  outlined in \Cref{subsec:natural-maps} (see \cite[Section 5.2]{RCDbary}), we associate to $\hat{f}$ a family of natural maps $\hat{F}_{s} \colon X \to \hat{Y}$ for $s > h(\overline{X})$. 
Explicitly, recall that $\hat{F}_s(x)$ is obtained by choosing a lift $\tilde{x} \in \overline{X}$, pushing forward the exponential measure $\mu_{\tilde{x}}^s = e^{-s \dist(\tilde{x}, \cdot)} \overline{\mathfrak{m}}$ to the $\CAT(0)$ universal cover $\til{Y}$ via the lifted map $\tilde{\hat{f}}$, then computing its unique $d^2$-barycenter in $\til{Y}$, and finally projecting this $\pi_1(\hat{Y})$-equivariant point down to $\hat{Y}$.
As $s \searrow h(\overline{X}) = h(\tilde{\hat{Y}})$, the maps $\hat{F}_{s}$ subconverge uniformly to a continuous limit map $\hat{F}_\infty \colon X \to \hat{Y}$. 
Because the Jacobian bounds of \Cref{prop:entropy-jacobian} imply that the asymptotic Lipschitz constant approaches $1$ as $s \searrow h(\til{Y})$, the limit map $\hat{F}_\infty$ is $1$-Lipschitz \cite[Lemma 5.1]{RCDbary}.
The manifold $\hat{Y}$ is negatively curved, so continuous maps into $\hat{Y}$ that are uniformly close are homotopic \cite[Lemma 5.5]{RCDbary}. 
Therefore, the limit map $\hat{F}_\infty$ is homotopic to $\hat{f}$. 
Then the homotopy invariance of the degree (\Cref{lem:homotopy-invariance}) implies that $|\deg(\hat{F}_\infty)| = |\deg(\hat{f})| = 1$. 

We now evaluate the mass of the push-forward of the orientation current $T$ under the $1$-Lipschitz map $\hat{F}_\infty$. 
Because $\hat{F}_\infty$ is $1$-Lipschitz, its metric Jacobian is bounded by $1$ almost everywhere. 
Combining \Cref{prop:degree-equivalence} with \Cref{lem:mass-pushforward}, we obtain:
\begin{align*}
    \Hn(\hat{Y}) &= |\deg(\hat{F}_\infty)|\, \Hn(\hat{Y}) = \bM\big(\hat{F}_{\infty\#} T\big) \\
    &\leq \int_X \Jac_x \hat{F}_\infty\, d\Hn(x) \\
    &\leq \int_X 1 \, d\Hn(x) = \Hn(X).
\end{align*}
The equality $\Hn(X) = \Hn(\hat{Y})$ implies that the inequalities above are strict equalities. 
Therefore $\Jac_x \hat{F}_\infty = 1$ for $\Hn$-almost every $x \in X$, meaning $\hat{F}_\infty$ is a volume-preserving map. 
Because $|\deg(\hat{F}_\infty)| = 1 \neq 0$, the map $\hat{F}_\infty$ is surjective by \Cref{cor:surjective}.

At this stage, we have a $1$-Lipschitz, surjective, volume-preserving map $\hat{F}_\infty \colon X \to \hat{Y}$ satisfying $\Hn(\hat{F}_\infty(X)) = \Hn(\hat{Y}) = \Hn(X)$ between two closed non-collapsed $\RCD$ spaces. 
By the Lipschitz-volume rigidity theorem for $\RCD$ spaces without boundary \cite[Theorem 1.8]{RCDbary} (extending the Alexandrov space rigidity result of Li--Wang \cite{li2014lipschitz}, and see also the concurrent work \cite[Corollary 1.4]{Li-Pasqualetto-Qian-2026}), we conclude that the map $\hat{F}_\infty$ is a metric isometry.

As $\hat{F}_\infty$ is a global isometry and $\hat{Y}$ is a smooth Riemannian manifold, $X$ inherits this smooth locally symmetric Riemannian structure. 
Moreover, the original continuous map $f$ is homotopic to the composition $F = p \circ \hat{F}_\infty \colon X \to Y$. 
Finally, because $\hat{F}_\infty$ is an isometry and $p$ is a Riemannian covering map of degree $k$, the map $F$ is a Riemannian covering map of degree $k = |\deg(f)|$. 
\end{proof}

\pagebreak

\subsection*{Report on Human and AI contributions}\label{app:AI}

I first proposed this problem in 2015, and initiated the collaboration developing the barycenter method for $\RCD$ spaces \cite{RCDbary}.
Recently, I formulated a degree theory suited to integral current spaces \cite{SuarezSerrato2026}, which, in combination with the recent work on orientability of metric spaces \cite{CMMP2026} contribute all the technical requirements for the main result here.
Knowing what was needed for this work, with the assistance of Gemini Pro 3.1 \& Argon, and Claude Opus 5.5 (in the sandboxed AI Commons resource at UCSB)  I completed the details of the current paper.
The AI's role was auxiliary, helping with \LaTeX \, formatting, style editing, argument review, and bibliographic verification.
The responsibility for the mathematical rigor and correctness is mine.

%

\providecommand{\bysame}{\leavevmode\hbox to3em{\hrulefill}\thinspace}

\end{document}